\documentclass[a4paper,11pt]{article}
\usepackage[utf8]{inputenc}

\usepackage{amssymb,amsmath,amsfonts}
\usepackage{graphicx}
\usepackage{color}
\usepackage{hyperref}
\usepackage{cleveref}
\usepackage{subfig}
\usepackage{bm}
\usepackage{listings}
\usepackage{subfig}
\newcount\colveccount
\newcommand*\colvec[1]{
        \global\colveccount#1
        \begin{pmatrix}
        \colvecnext
}
\def\colvecnext#1{
        #1
        \global\advance\colveccount-1
        \ifnum\colveccount>0
                \\
                \expandafter\colvecnext
        \else
                \end{pmatrix}
        \fi
}

\newcommand{\ymp}[1]{{\color{black}{#1}}}
\newcommand{\ympp}[1]{{\color{black}{#1}}}
\newlength{\figwidth}
\figwidth=\linewidth

\newcommand{\qed}{\nobreak \ifvmode \relax \else
      \ifdim\lastskip<1.5em \hskip-\lastskip
      \hskip1.5em plus0em minus0.5em \fi \nobreak
      \vrule height0.75em width0.5em depth0.25em\fi}

\usepackage{cite} % Make references as [1-4], not [1,2,3,4]
\usepackage{url}  % Formatting web addresses  
\usepackage{ifthen}  % Conditional 
\usepackage{multicol}   %Columns
\usepackage[utf8]{inputenc} %unicode support

\title{Weakly Coupled Oscillators in a Slowly Varying World}
\author{Youngmin Park\footnote{Corresponding author. Email yop6@pitt.edu} \& Bard Ermentrout \\ Department of Mathematics\\University of Pittsburgh\\ Pittsburgh PA 15260}

\begin{document}

\maketitle
%% \tableofcontents
\begin{abstract}
We extend the theory of weakly coupled oscillators to incorporate slowly varying inputs and parameters. We employ a combination of regular perturbation and an adiabatic approximation to derive equations for the phase-difference between a pair of oscillators.  We apply this to the simple Hopf oscillator and then to a biophysical model. The latter represents the behavior of a neuron that is subject to slow modulation of a muscarinic current such as would occur during transient attention through cholinergic activation.  Our method extends and simplifies the recent work of Kurebayashi \cite{kureb13} to include coupling. We apply the method to an all-to-all network and show that there is a waxing and waning of synchrony of modulated neurons. 

\end{abstract}

\bigskip

\section{Introduction}
The theory of weakly coupled oscillators \cite{kuramoto84,ermentrout1981,ek84} has served very well as a predictor of the dynamics in networks of coupled neural oscillators  (for a comprehensive review, see \cite{schwemmer}).  In the application of this theory, one generally assumes that, while the oscillators may have different intrinsic frequencies, these frequencies are fixed as are the uncoupled limit-cycle oscillators. However, more generally, local regions of the nervous system are constantly modulated by extrinsic inputs and by slow processes such as the accumulation of extracellular ions.  Thus, synchronization and other properties are likely to change due to this modulation which can change the frequency, conductances, and even the synapses within an oscillatory network \cite{pfeuty}.  

Neuronal modulators such as acetylcholine, norepinephrine \cite{mccormick1992}, and dopamine \cite{gorelova} are known to alter the firing properties of neurons.  These properties, in turn, could alter the synchronization behavior of neurons and more formally, the form that the weak coupling equations take. One of the key components to understanding synchronization of neuronal oscillators is the phase response curve (PRC) which describes how the phase of an oscillator is shifted by the timing of inputs. The PRC plays the key role in determining whether or not a pair of coupled neuronal oscillators will synchronize or not.  In \cite{stiefel1,stiefel2}, the authors directly demonstrated that cholinergic modulation of a cortical pyramidal neuron had a profound effect on the shape of the PRC.  Acetylcholine is known to directly act on the so-called M-type potassium current and in \cite{pascal}, they showed how changing the strength of this current made a huge difference on the shape of the PRC as well as in the ability of synaptically coupled neurons to synchronize. Neuronal properties are also affected by the extracellular milieu, notably, concentration of extracellular potassium which can profoundly alter excitability of neurons \cite{cressman}. Rubin et al \cite{rubinrubin} showed that the synchronization between two coupled neurons was strongly dependent on the mean concentration of extracellular potassium.  Jeong and Gutkin \cite{jeong07} showed the changes in the reversal potential of GABAergic conductances changed the ability of neurons to synchronize; the reversal potential is mainly driven by extracellular chloride.  Thus, since many neuromodulators as well as the ionic concentrations are constantly changing, it is important to see how this time varying environment alters the ability of neurons to synchronize.

In two recent papers Kurebayashi et al \cite{kureb13,kureb15}, extended the notion of phase reduction to oscillators that are subject to large slowly varying parameters. They demonstrated that the evolution of the phase depended, not just on the instantaneous frequency of the oscillator, but, also on the rate of change of the slowly varying parameter.  In this paper, we re-derive the phase equation in \cite{kureb13} by using the method of adiabatic invariance \cite{keener} (Chapter 12.1.2) and incorporate the slow variation of parameters into weak coupling of oscillators using the Fredholm alternative.  Thus we have a theory to predict synchrony and antiphase along with stability, in the presence of a slowly varying parameter. Moreover, because we only assume the parameter to be slowly varying, our theory is shown to accurately predict phase differences with periodic, quasi-periodic, and stochastic slowly varying parameters.  

We first derive the equations for the phases and the phase-differences for a pair of coupled oscillators that are subject to slow changes in a parameter.  Next, we apply the theory to the Hopf oscillator (so-called $\lambda-\omega$ system, \cite{kopellhoward}) where all of the required functions for our analysis can be exactly derived.  We then consider a biophysical Hodgkin-Huxley model for pyramidal neurons (the simplified ``Traub'' model \cite{pascal}). This model includes an M-type potassium current, so in our analysis and simulations, we allow the conductance to slowly change as a model for cholinergic modulation.  We conclude with a discussion and contrast the results with fast modulation.

\section{Methods}
\subsection{Weakly Coupled Oscillators With a Slowly Varying Parameter}
We consider a pair of weakly coupled slowly-varying oscillators:
\begin{eqnarray}
\label{eq:coupled}
\frac{dX^a}{dt} &=&   F(X^a,q\,\!(\epsilon t)) + \epsilon G_a(X^b,X^a) \\
\frac{dX^b}{dt} &=&   F(X^b,q\,\!(\epsilon t)) + \epsilon G_b(X^a,X^b) \nonumber
\end{eqnarray}
where $0<\epsilon\ll 1$ is a small parameter.  We assume that the slowly varying parameter, $q\,\!$, lies in an interval $Q:=[q\,\!^-,q\,\!^+]$ such that for each $q\,\!\in Q$, the system:
\[
\frac{dX}{dt} = F(X,q\,\!)
\]
has an asymptotically stable limit cycle with frequency $\omega(q\,\!)$. The period of the oscillators is just $T(q\,\!)=2\pi/\omega(q\,\!).$    Thus, each of the two oscillators is modulated by a common slowly varying signal, $q\,\!(\epsilon t)$ that can alter the shape and frequency of the rhythm but does not destroy its existence.   The functions $G_{a,b}$ represent the weak coupling between the two oscillators.  If there is no modulation of the oscillations, then, we can regard equation (\ref{eq:coupled}) as a standard weakly coupled system. However, the slow modulation changes the dynamics and interactions in a way that we will now demonstrate.  We point out that \cite{kureb13} derived the phase modulation for a single slowly varying oscillator using successive changes of variables and showed that the ``naive'' phase approximation was not valid. \ymp{More precisely, the term $\beta(\tau)$ defined in equation \eqref{eq:beta} accounts for possibly large variations in the slowly varying parameter. Therefore, omitting this term \ymp{(the ``naive'' approximation)} from equations \eqref{eq:tha}--\eqref{eq:thb} results in a poor phase approximation in the case of a single forced oscillator. Here, we introduce a simpler way to derive the same equations using a standard adiabatic approximation \cite{keener}, and extend the results to coupled oscillators. In the coupled oscillator case, the ``naive'' phase approximation is valid.} 

Clearly, there are two time scales in this problem, a slow time scale, $\tau = \epsilon t$ and a fast time scale, $s$ that is related to $t$.  To put everything on a similar fast time scale, we generally allow that
\[
\frac{ds}{dt} = g(\tau,\epsilon)
\]
and expand to get a relationship between $s$ and $t$ that is $\tau-$dependent.  However, we need only terms in the lowest order in the fast time scale, so we will simply cut to the chase and write $s=\omega(q\,\!)t$, so that the oscillators are all $2\pi-$periodic in $s$. (For the time being, we have suppressed the implicit $\tau$ dependence of the fast time scale, by just putting the parameter $q\,\!$ in the frequency, but, in fact, $q\,\!$ is just shorthand notation for $q\,\!(\tau).$)

Before continuing with the perturbation, we introduce some additional notation.  Let $U_0(s,q\,\!)$ be the limit cycle solution to the uncoupled system
\begin{equation}
\label{eq:uncoupled}
 \omega(q\,\!) \partial U/\partial s=F(U,q\,\!),
\end{equation}
and let $A(s,q\,\!):=D_UF(U_0,q\,\!)$ be the linearization of the uncoupled system evaluated at the limit cycle. \ymp{By taking the derivative with respect to $s$ on both sides of equation \eqref{eq:uncoupled}, we see that} the linear equation
\begin{equation}
\label{eq:lin}
L(s,q\,\!) Y:= \omega(q\,\!) \frac{\partial Y}{\partial s} - A(s,q\,\!)Y =0
\end{equation}
has a periodic solution given by $\partial_s U_0(s,q\,\!)$ where the notation, $\partial_s$ means differentiation with respect to the first component in $U_0$.  Associated with the set of $2\pi-$periodic functions is an inner product defined as 
\[
\langle Y_1(s),Y_2(s) \rangle = \int_0^{2\pi} Y_1(s)\cdot Y_2(s)\ ds,
\]
where $Y_1\cdot Y_2$ is the standard Euclidean dot product.  With this inner product,  the linearized equation has a well-defined adjoint operator:
\ymp{\begin{equation}
 L^*(s,q\,\!) Y:= \omega(q\,\!) \frac{\partial Y}{\partial s} + A^T(s,q\,\!)Y ,
\end{equation}
from which we attain the adjoint equation,}
\[
\omega(q\,\!) \frac{\partial Y^*}{\partial s} = -A^T(s,q\,\!)Y^*
\]
where $A^T$ is the transpose.  There is a unique $2\pi-$periodic solution to the \ymp{adjoint equation}, which we call $Z(s,q\,\!)$ such that
\[
Z(s,q\,\!)\cdot \frac{\partial U_0(s,q\,\!)}{\partial s} =1.
\]
Finally, the linearization has the Fredholm alternative property \ymp{\cite{keener}}. That is, there is a $2\pi-$periodic solution to:
\begin{equation}
\label{eq:inhom}
L(s,q\,\!)Y = b(s)
\end{equation}
where $b(s)$ is $2\pi-$periodic and $L$ is defined in (\ref{eq:lin}) if and only if
\[
\int_0^{2\pi} Z(s,q\,\!)\cdot b(s)\ ds = 0.
\]
With these preliminaries defined, we are now ready to analyze weak coupling of slowly varying oscillators.  

We assume that the solutions to equation (\ref{eq:coupled}) can be expressed in a series in $\epsilon$ and have the form 
\begin{eqnarray*}
X^a(t,\epsilon) &=&  X^a_0(s,\tau) + \epsilon X^a_1(s,\tau) + \ldots \\
X^b(t,\epsilon) &=&  X^b_0(s,\tau) + \epsilon X^b_1(s,\tau) + \ldots.
\end{eqnarray*}
To lowest order, we must have
\[
\omega(q\,\!(\tau)) \frac{\partial X^{a,b}_0(s,\tau)}{\partial s} = F(X^{a,b}_0(s,\tau),q\,\!(\tau)). 
\]  
The $2\pi-$periodic limit cycle solution to this problem is
\[
X^{a,b}_0(s,\tau) = U_0(s+\theta^{a,b}(\tau),q\,\!(\tau))
\]
where $\theta^{a,b}(\tau)$ are slowly varying arbitrary phase shifts due to the \ymp{time-translation invariance} of the limit cycle.   Our goal is to now derive equations for the slow evolution of the phase.  Using the chain rule, we see that  $d/dt = \omega \partial /\partial s + \epsilon \partial/\partial \tau.$ Thus, after a bit of rearranging, the next order equations are:
\begin{eqnarray}
\label{eq:ord1}
L(s,q(\tau)) X^a_1(s,\tau) &=&
-\partial_s U_0(s+\theta^a(\tau),q\,\!(\tau))\frac{\partial \theta^a}{\partial \tau} \\
{} &{}& -\partial_q\,\! U_0(s+\theta^a(\tau),q\,\!(\tau))\frac{dq\,\!}{\partial \tau} \nonumber \\
{} &{}& +G_a[U_0(s+\theta^b(\tau),q\,\!(\tau)),U_0(s+\theta^a(\tau),q\,\!(\tau))], \nonumber
\end{eqnarray}
where $L(s,q\,\!(\tau))$ is defined by (\ref{eq:lin}). There is a similar equation for  $X^b_1(s,\tau).$ Finally, we see that this equation has the form of (\ref{eq:inhom}), so that there is a $2\pi-$periodic solution if and only if the right-hand side is orthogonal to the adjoint. This leads to the following equations for the phases:
\begin{eqnarray}
\label{eq:tha}
\frac{\partial \theta^a}{\partial \tau} &=& -\beta(\tau) + h_a(\theta^b-\theta^a,\tau) \\ 
\label{eq:thb}
\frac{\partial \theta^b}{\partial \tau} &=& -\beta(\tau) + h_b(\theta^a-\theta^b,\tau), 
\end{eqnarray}
where,
\begin{eqnarray}
\label{eq:beta}
\beta(\tau) &=& \int_0^{2\pi} \ymp{Z}(s,q\,\!(\tau))\cdot\partial_q\,\! U_0(s,q\,\!(\tau))\frac{\partial q\,\!}{\partial \tau} \ ds \\
\label{eq:H}
h_{a,b}(\phi,\tau) &=& \int_0^{2\pi} \ymp{Z}(s,q\,\!(\tau))\cdot G_{a,b}[U_0(s+\phi,q\,\!(\tau)),U_0(s,q\,\!(\tau))] \ ds.
\end{eqnarray}
The extra $\beta(\tau)$ term arises due to the fact that the parameter $q\,\!$ is slowly varying. Notice in its definition through equation (\ref{eq:beta}), that it is proportional to the time derivative of $q\,\!(\tau).$ This is the term that \cite{kureb13} emphasized in their analysis. That is, we recover their results if we ignore coupling. \ymp{We remark that the phase-interaction functions, $h_{a,b}(\cdot)$ are exactly those that would be obtained from standard weak-coupling theory with all parameters held fixed.} In absence of coupling, the total phase evolves as
\[
\theta(t) = \theta(0) + \omega(q\,\!(\tau)) t - \epsilon \int_0^t \beta(\epsilon t')\ dt'.
\]
If both oscillators are subject to the exact same slowly varying inputs, then the $\beta$ term becomes irrelevant to their phase difference, $\phi:=\theta^b-\theta^a$ which satisfies the simple scalar slowly varying equation \ymp{\cite{schwemmer}}:
\begin{equation}
\label{eq:phi}
\frac{d \phi}{d\tau} = h_b(-\phi,\tau)-h_a(\phi,\tau):=G(\phi,\tau)
\end{equation} 
Equation~(\ref{eq:phi}) will be our main tool for comparing the phase reduced model to the full model. We remark that the interaction functions are $\tau-$dependent, so that the reduced system for the phase-difference is no longer autonomous and we will not be able to write exact solutions.  However, if $h_a=h_b$, then the right-hand side of equation (\ref{eq:phi}) is, for each $\tau$, an odd periodic function of $\phi$, so that $\phi=0,\pi$ will always be equilibrium points. That is $G(0,\tau)=G(\pi,\tau)=0$ for all $\tau.$ In this symmetric case, we define $H_{odd}$ to be the odd part of the function $h_a$.

\subsection{Mode Truncation}
In order to study the phase-reduced equations, we need to get formulae for the $\tau-$dependent interaction functions, $h_{a,b}$.  For our first application of the method, these functions are explicitly computable since the oscillation and the adjoint are simple sine and cosine functions.  However, for the neural model that we also study (and which gives more interesting results), we need to somehow approximate the required slowly varying functions. To this end, we use XPPAUT \cite{xpp}  to numerically compute the adjoint and interaction functions. For each interaction function we perform a mode truncation (that is, we keep just a few of the Fourier terms). We finish the approximation by deriving a $q$-dependent equation for the coefficients of the Fourier series expansion; generally piecewise-linear.  
This approximation serves particularly well for our problem because only two sine coefficients are required to preserve the change in synchrony, bistability between synchrony and antiphase ($\pi-$phase difference), and other interesting phenomena.

\section{Results}
We apply our theory to the Lambda-Omega system and a modified Traub model with adaptation. For each model, we consider three types of slowly varying parameters, which we briefly discuss before delving into the details of each model. \ympp{We remark that all figure code and relevant data files are available on github at \url{https://github.com/youngmp/park_and_ermentrout_2016}}

\subsection{Slowly Varying Parameters}
The slowly varying parameter, $q(\tau)$, is explicitly written as three types of slowly varying parameters: periodic ($q_p$), quasi-periodic ($q_{qp}$), and stochastic ($q_s$):
\begin{equation}
\label{eq:slowv}
\begin{split}
 q_p(\tau) &:= q_0 + q_1 \cos(f\tau),\\
 q_{qp}(\tau) &:= q_0 + (q_1/2)(\cos(f\tau)+\cos(f\tau\sqrt{2})),\\
 q_s(\tau) &:= q_0 + q_1 z(\tau).
\end{split}
\end{equation}
The terms $q_i$, $f$, and $\varepsilon$ depend on the system. \ymp{For the Lambda-Omega system, we choose $f=1$ and various combinations of $q_0$ and $q_1$ because the choice of $q_i$ affects the \ympp{asymptotic} dynamics. Surveying multiple values of $q_i$ provides a more complete demonstration of the dynamics and the accuracy of our theory. For the Traub model, we chose by default $q_0=0.3$, $q_1=0.2$, and $f=5$ unless otherwise stated (as in figure \ref{fig:trb_p_n}). The default choice of parameters represents a biophysically realistic parameter range, while the slightly different parameter choice in figure \ref{fig:trb_p_n} demonstrates the accuracy of our theory when we avoid slow stability changes.}

The noisy parameter, $z$, is an Ornstein-Uhlenbeck (OU) process satisfying the stochastic differential equation
\begin{equation*}\label{eq:ou}
 \mu \mathrm{d}z=-z\mathrm{d}t+\sqrt{\mu}\mathrm{d}W,
\end{equation*}
where $\mu=1000$. The raw random noise data is normalized so that
\begin{equation*}
 z(\tau) \in [-1,1], \quad \forall \tau,
\end{equation*}
and this data is used in all noisy simulations. The OU data may be reproduced by using XPP seeds 1--4. 

\subsection{Lambda-Omega System}
We first apply our result to the $\lambda-\omega$ system \cite{kopellhoward} with weak diffusive coupling,

\begin{equation}\label{eq:lamom2}
 \left( \begin{matrix}\dot x_j \\ \dot y_j\end{matrix} \right) = \left( \begin{matrix}\lambda(r_j) & -\omega(r_j,q(\tau)) \\ \omega(r_j,q(\tau)) & \lambda(r_j)\end{matrix} \right) \left( \begin{matrix} x_j \\ y_j\end{matrix} \right) + \varepsilon \left( \begin{matrix} 1 & -\ymp{\kappa} \\ \ymp{\kappa} & 1  \end{matrix} \right)\left( \begin{matrix} x_k-x_j \\ y_k-y_j  \end{matrix} \right)
\end{equation}
where $j,k=1,2$ and $k\ne j$;  $r_j:=\sqrt{x_j^2+y_j^2}$; and
\begin{eqnarray*}
\lambda(r) &=& 1-r^2 \\
\omega(r,q) &=& \ymp{1 + q(r^2-1)}.
\end{eqnarray*}

\ymp{When $\varepsilon=0$, equation \eqref{eq:lamom2} is equivalent to the Hopf oscillator in polar coordinates,
\begin{eqnarray*}
\dot r_j &=& r_j\lambda(r_j)\\
\dot \theta_j &=& \omega(r_j,q(\tau)).
\end{eqnarray*}}

One can verify the limit cycle for uncoupled system is
\begin{equation*}
 U_0(s,\tau) = [\cos(s),\sin(s)]^T,
\end{equation*}
and the solution to the adjoint equation (the iPRC) is
\[
 Z(t) = [q(\tau)\cos(t) - \sin(t),q(\tau)\sin(t)+\cos(t)]^T.
\]
 Finally, Equation~(\ref{eq:phi}) for the $\lambda-\omega$ system is
\begin{equation}\label{eq:lamom_hodd}
 \frac{d\phi}{d\tau} = 2\left(\ymp{\kappa} q(\tau)- 1\right)\sin(\phi).
\end{equation}
Note that synchrony \ymp{($\phi=0$)} is indeed a fixed point of Eq.~\ref{eq:lamom_hodd}. For a brief stability analysis, we note that equation \eqref{eq:lamom_hodd} is a separable equation and solve for an implicit solution the differential equation:
\begin{equation}
 \tan(\phi/2) = c \exp\left[\int_0^\tau (\ymp{\kappa} q(s) - 1) ds\right].
\end{equation}
\ymp{We can write the inside of the exponential as 
$\tau Q(\tau)$ where $Q(\tau)=(1/\tau)\int_0^\tau [\kappa q(s)-1]\ ds.$  $Q(\tau)$ is the running average of the integrand. Since the integrand is bounded and continuous, the limit of $Q(\tau)$ exists as $\tau\to\infty.$  If this limit is positive then the exponential diverges to $+\infty$ and the phase $\phi\to\pm\pi.$  Similarly, if the limit is negative, then the exponential goes to zero and $\phi$ converges to 0 as well.}

\begin{figure}[h!]
\centering
\includegraphics[width=1\textwidth]{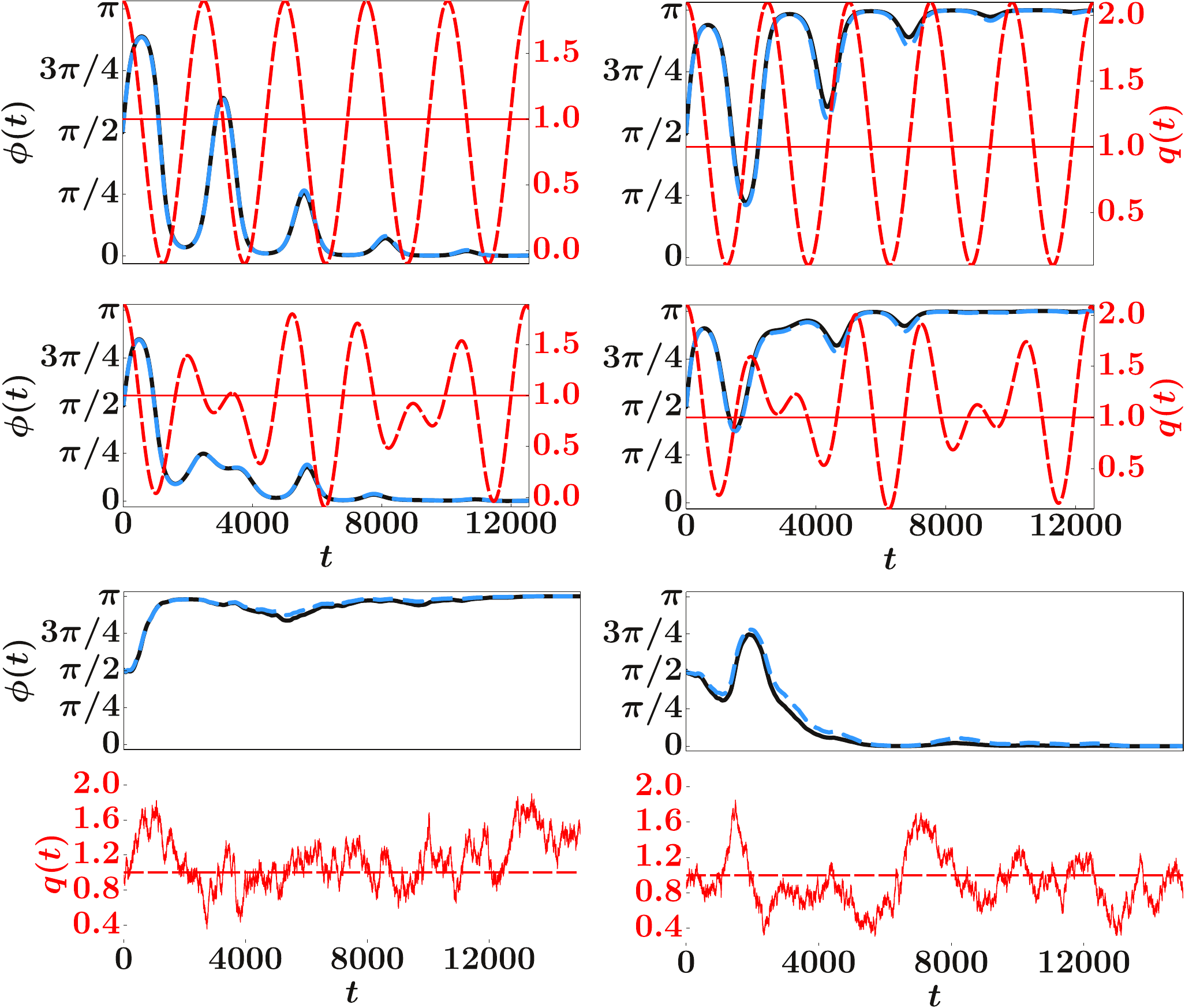}
 \caption{Periodic ((a),(b)), quasi-periodic ((c),(d)) and stochastic ((e),(f)) slowly varying parameters. The phase difference $\phi = \theta_2 - \theta_1$ theory (light blue dashed line) is plotted on top of experiment (black solid line). The slow periodic parameter is shown as a dashed red line. The horizontal line represents the parameter value $q$ at which there is onset or offset of synchrony. For all subfigures, $\varepsilon=0.0025$, $f=\ymp{\kappa}=q_1=1$. (a) Periodic, $q_0=0.9$, $\mathbb{E}[q]<1$. (b) Periodic, $q_0=1.1$, $\mathbb{E}[q]>1$. (c) Quasi-periodic, $q_0=0.9$, $\mathbb{E}[q]<1$. (d) Quasi-periodic, $q_0=1.1$, $\mathbb{E}[q]>1$. (e) Stochastic (OU), $q_0=0.85$, $\mathbb{E}[q]=0.913$, XPP seed 2. (f) Stochastic (OU), $q_0=0.9$, $\mathbb{E}[q]=1.145$, XPP seed 1.}\label{fig:lamom}
\end{figure}

Figure \ref{fig:lamom} shows the result of simulating equation (\ref{eq:lamom2}) for different functions of $q(\tau)$. In the left column (labeled a,c,e) the mean value of $q(\tau)$ is less than 1, so that we expect that the phase differences will go to synchrony. In the right panels, the mean value of $q(\tau)>1$ so that the theory predicts that  phase-differences will go to $\pi.$  This is clearly evident from the simulations of the full model. Furthermore, the approach to equilibrium predicted by the phase model is almost identical to that of the full simulations. There is very little error even in the stochastic cases (panels e,f). Even though this is a highly nonlinear system, the system goes to the stable state that is appropriate for the {\em average} of the slowly varying parameter.  If we break the homogeneity, then the dynamics is more complex and interesting.

\subsubsection*{Heterogeneities.} 
\begin{figure}
\centering
\includegraphics[width=\textwidth]{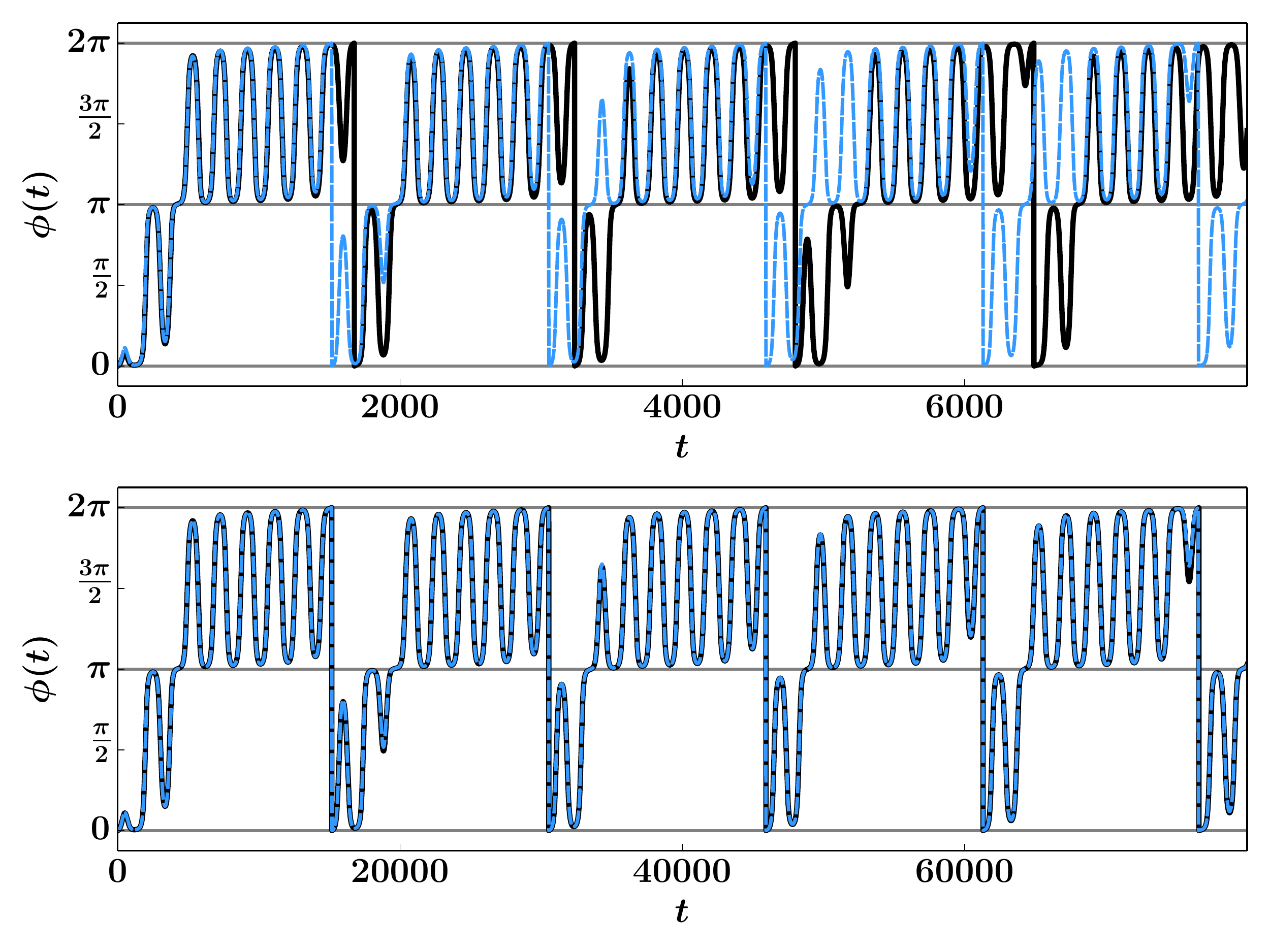}
\caption{The effects of inhomogeneity on slowly modulated solutions. The slowly varying parameter is chosen to be periodic with $q_0=1.1,q_1=2,f=1.3, d=0.05$, $\ymp{\kappa}=1$.  (a) $\epsilon=0.025$ (b) $\epsilon=0.0025$. Black is the full model and blue dashed is the phase-reduced model. \ympp{The solid gray lines at $\phi(t)=0\equiv 2\pi$ ($\phi(t) = \pi$) represent synchrony (anti-phase).}}
\label{fig:loinh}
\end{figure}

\ymp{So far, the model derivation assumes that both oscillators are identical. However, in general, this will not necessarily be the case. If the differences are $O(\epsilon)$ (that is, small) then we will get some additional terms in the phase equations.  To account for this, in general, we add terms to equation (\ref{eq:coupled}) of the form:
\begin{equation}
\label{eq:extra}
\epsilon f_{a,b}(X^{a,b},\tau)
\end{equation}
where we could also include some $\tau$-dependence in the heterogeneity. For example, in the $\lambda-\omega$ system, we could set
\[
\omega_{a,b}(r,q) = 1 + q(\tau)(r^2-1) + \epsilon [d_{a,b} + c_{a,b}(\tau)(r^2-1)]
\]
where the subscripts refer again to the two oscillators. Here, the parameters $d_{a,b}$ are just constants that affect the baseline frequency and $c_{a,b}(\tau)$ are modulatory.  With the addition of the terms (\ref{eq:extra}), the phase equations we get are like equations (\ref{eq:tha}-\ref{eq:thb}) but have additional terms:
\begin{eqnarray*}
\frac{\partial \theta^a}{\partial \tau} &=& -\beta(\tau) + \eta_a(\tau) + h_a(\theta^b-\theta^a,\tau) \\ 
\frac{\partial \theta^b}{\partial \tau} &=& -\beta(\tau) + \eta_b(\tau)+ h_b(\theta^a-\theta^b,\tau), 
\end{eqnarray*}
where
\[
\eta_{a,b}(\tau) = \int_0^{2\pi} Z(s,q(\tau))\cdot f_{a,b}(U_0(s,q(\tau)),\tau)\ ds.
\]
Subtracting the two equations yields the more general phase equation with heterogeneities:
\begin{equation}
\label{eq:phihet}
\frac{d \phi}{d\tau} = \eta_b(\tau)-\eta_a(\tau)+h_b(-\phi,\tau)-h_a(\phi,\tau):=G(\phi,\tau).
\end{equation} 
We have still eliminated the commmon $O(1)$ slow variation $\beta(\tau)$, but the explicit heterogeneities appear threough the differences $\eta_b(\tau)-\eta_a(\tau).$ 
 
With this extension,} we now alter the simple model  by introducing a small frequency difference in the oscillators. For oscillator 2, we replace $\omega(r,q)=1+q(1-r^2)$ with  $\omega(r,q)=1+\epsilon d + q(1-r^2)$, so that in absence of coupling, there is an order $\epsilon$ frequency difference, $\epsilon d.$ In this case the equation for $\phi$ becomes
\begin{equation}
\label{eq:loinh}
 \frac{d\phi}{d\tau} = d+2\left(\ymp{\kappa} q(\tau)- 1\right)\sin(\phi).
\end{equation}
This means that $\phi(\tau)$ will no longer generally approach a steady state.  In figure \ref{fig:loinh} we show two simulations with different values of $\epsilon$ when there is a slight difference in frequency.  For $\epsilon=0.025$, the solutions match for most of the time, but there are places in each segment, where the solutions are about $\pi$ out of phase. On the other hand, when we reduce $\epsilon$ by  factor of 10, the solutions to the full model and the phase model are indistinguishable.

\begin{figure}
\centering
\includegraphics[width=\textwidth]{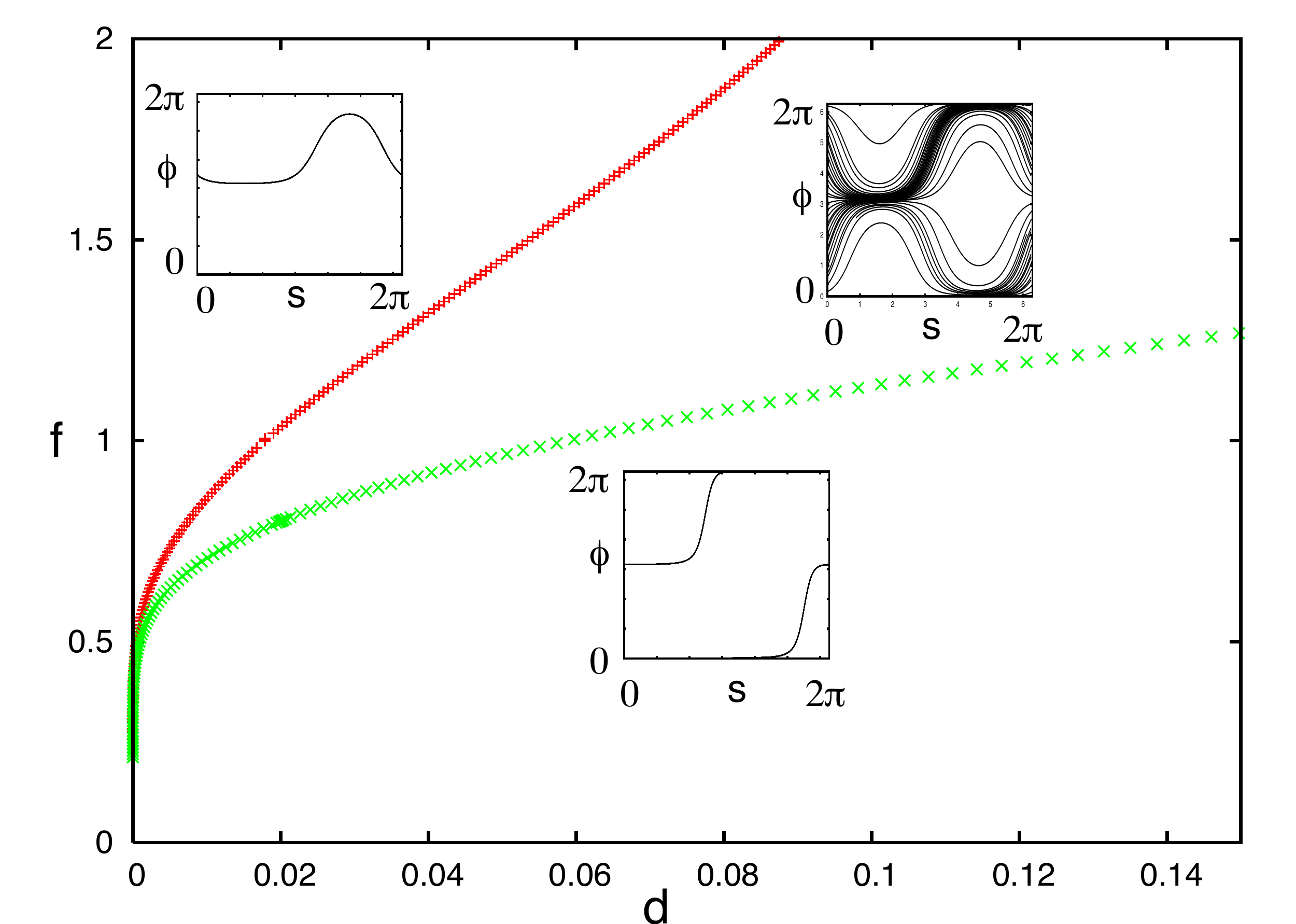}
\caption{The behavior of equation (\ref{eq:loinh}) for periodic modulation as a function of the homogeneity, $d$ and the modulation frequency, $f$. Green points show the border for $1:1$ locking; red points show the border for $0:1$ locking and between these are mixed solutions. Typical phase-planes are shown in each region.}
\label{fig:lown}
\end{figure}

\ymp{As the effects of heterogeneities are rather interesting, even in this simple case, we will now examine equation (\ref{eq:loinh}) in more detail in order to explain the behavior in figure \ref{fig:loinh}.  We can rewrite equation (\ref{eq:loinh}) as a system with the time rescaled:
\begin{eqnarray*}
\phi' &=& (d+2(q(s)-1)\sin\phi)/f \\
s' &=& 1.
\end{eqnarray*} 
This is an equation on the torus and so the behavior is fairly restricted; in particular, the ratio $\rho=\lim_{\tau\to\infty}\phi/s$, called the rotation number is a continuous function of the parameters. We find three different behaviors as the inhomogeneity $d$ and the frequency $f$ vary.  Figure \ref{fig:lown} shows the behavior as these parameters are varied. In the upper left part of the diagram (high frequency), above the red curve, $\phi(\tau)$ has a winding number of 0. That is, $\phi(\tau+2\pi)=\phi(\tau).$  This means that the phase-difference, $\phi$ between the two oscillators is bounded between two values and one oscillator is consistently ahead of the other.  In the lower right part of the diagram (low frequency), $\phi(\tau+2\pi)=\phi(\tau)+2\pi$, that is, $\phi$ has winding number 1.  This means that the phase-difference between the two neurons stays close to 0 for about half a cycle and close to $\pi$ for the other half and makes these switches rapidly and periodically; it does not get ``stuck'' at synchrony or anti-phase.  Finally, the middle region (and also the choice used in figure \ref{fig:loinh}) shows that the the phase makes rapid transition, first between $\pi$ and $2\pi$ and then between $\pi$ and 0.  This explains the switching back and forth observed in figure \ref{fig:loinh}.  In sum, heterogeneity (even in the simplest form) can add good deal of complexity to the dynamics.}

\subsection{Traub Model with Adaptation}
The membrane potential dynamics of the Traub model, $V$, satisfies
\begin{equation}\label{eq:trb_v}
 \begin{split}
  C\dot V &= -g_{Na}m^3h(V-E_{Na}) - (g_k n^4+q(\tau) w)(V-E_k) - g_l(V-E_l)+I\\
  &\equiv f(V, q(\tau)),
 \end{split}
\end{equation}
where $q(\tau)$ is the slowly varying parameter with $q \in [0.1,0.5]$, $q_0 = 0.3$, $q_1 = 0.2$, and gating variables $n, m, h, w$ satisfying
\begin{equation*}\label{eq:trb_gates}
\begin{split}
 \dot n &= a_n(V)(1-n) - b_n(V)n,\\
 \dot m &= a_m(V)(1-m) - b_m(V)m,\\
 \dot h &= a_h(V)(1-h) - b_h(V)h,\\
 \dot w &= (w_{\infty}(V)-w)/t_w(V).
\end{split}
\end{equation*}
We introduce weak coupling by adding a synaptic conductance
\begin{equation*}\label{eq:trb_coupling}
\begin{split}
 \frac{dV_1}{dt} = f(V_1,q(\tau)) +\varepsilon g s_2(E_{syn} - V_1),\\
 \frac{dV_2}{dt} = f(V_2,q(\tau)) +\varepsilon g s_1(E_{syn} - V_2),
 \end{split}
\end{equation*}
where $s_i$ is the synaptic conductance of $V_i$ and satisfies
\begin{equation*}
 \dot s_i = \alpha(V_i)(1-s_i)-s_i/\tau_s.
\end{equation*}

\ymp{Here $\alpha(V)=4/(1+\exp(-v/5))$ \cite{wb96}}. Adaptation in this model is controlled by the magnitude of the $M-$type potassium current. This low-threshold, slow current can drastically affect the dynamics of the Traub model \cite{pascal} changing it from Class I excitable (oscillation arises via a saddle-node infinite cycle or SNIC)  to class II excitable (oscillation arises via a sub-critical Hopf bifurcation).  Because of that change, the adjoint, $Z(t)$ can also drastically change \cite{brown,ebn} and thus, the interaction function will also be strongly affected.   Biologically, this current is quite important since it is altered by acetylcholine, a neuromodulator.  Thus, since neuromodulators tend to operate at much slower time scales than the firing rates of neurons, the slow alteration of the $M-$type potassium current is an ideal example of the methods we have developed in this paper.

\begin{figure}
\centering
\includegraphics[width=\textwidth]{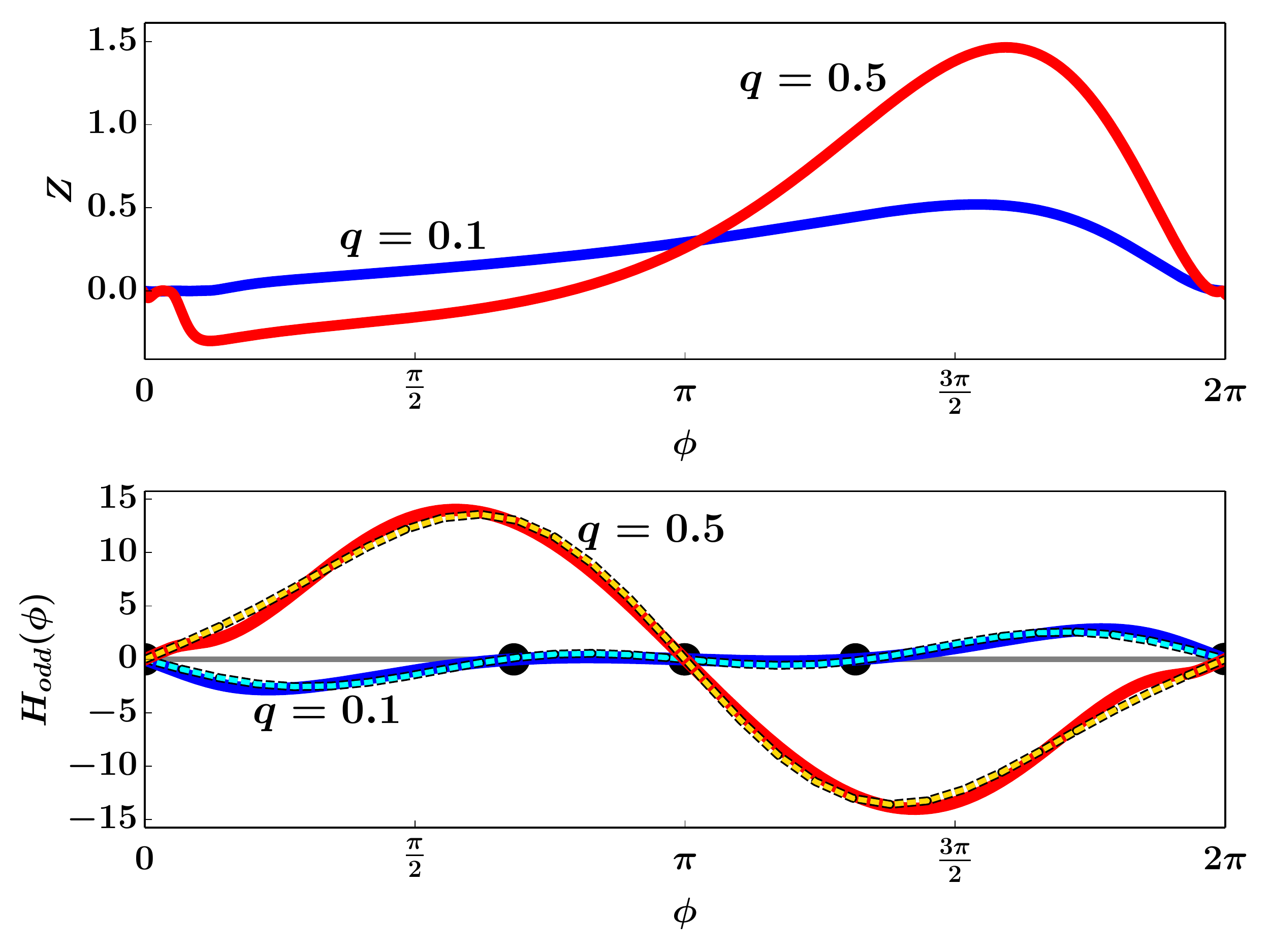}
\caption{The Traub model (equation (\ref{eq:trb_v}) for two fixed values of $q$, the M-type potassium current.  (a) the adjoints for $q=0.1$ and $q=0.5$; (b) The odd part of the interaction functions for $q=0.1,0.5$ (thick lines) and two-term sine fit (thin dashed lines) of $H_{odd}(\phi)$. \ymp{Zero crossings in (b) are denoted by black circles.}}

\label{fig:trbfg}
\end{figure}

Figure \ref{fig:trbfg} shows the results of a numerical computation of the adjoint and the odd part of the interaction function.  For small values of the M-current ($q=0.1$) the adjoint (a) is almost strictly positive which is typical for so-called Class I excitable systems where the periodic orbit arises as a SNIC.  On the lower panel (b) we see that the $H_{odd}(\phi)$ is small and that synchrony is unstable. (Recall that the phase model satisfies, $\phi'=-2H_{odd}(\phi)$, so that a negative (positive) slope at an equilibrium is unstable (stable).)   Anti-phase (\ymp{$\phi=\pi$}) is also unstable, but there are two stable fixed points that are near anti-phase.   When there is sufficient M-current ($q=0.5$), the adjoint has a large negative lobe right after the spike. This qualitative change in the shape of the adjoint leads to the stabilization of the synchronous state (panel b).  Thus, as $q$ is varied from a low to high value, we expect that the phase-difference will move toward synchrony (at high values) and away from synchrony (at lower values). Panel b also shows that a two-term sine approximation is reasonable and captures the qualitative (and to some extent, quantitative) shape of the functions. In particular, the full $H_{odd}(\phi)$ and the two-term sine approximation have the same equilibrium point properties.  For this reason, we make a simple linear interpolation using a two term sine expansion of the interaction as $q$ slowly varies. The approximation is thus:
\begin{equation*}
 -2H_{\text{odd}} \approx 2(b_1(q(\tau))\sin(\phi) + b_2(q(\tau))\sin(2\phi)),
\end{equation*}
\ymp{where a linear approximation to $b_i(q)$ passing through the points $(0.1,b_i(0.1))$ and $(0.3,b_i(0.3))$} predicts onset and offset of synchrony sufficiently well:
\begin{equation*}
 b_{i}(q) = 5(\hat b_i(0.3) - \hat b_i(0.1))q + 1.5 \hat b_i(0.1) - 0.5 \hat b_i(0.3),\quad i=1,2.
\end{equation*}
The number represented by $\hat b_i(x)$ is the actual coefficient value at $q=x$ (see \cref{app:fourier_coeff}).

\begin{figure}
\centering
\includegraphics[width=\textwidth]{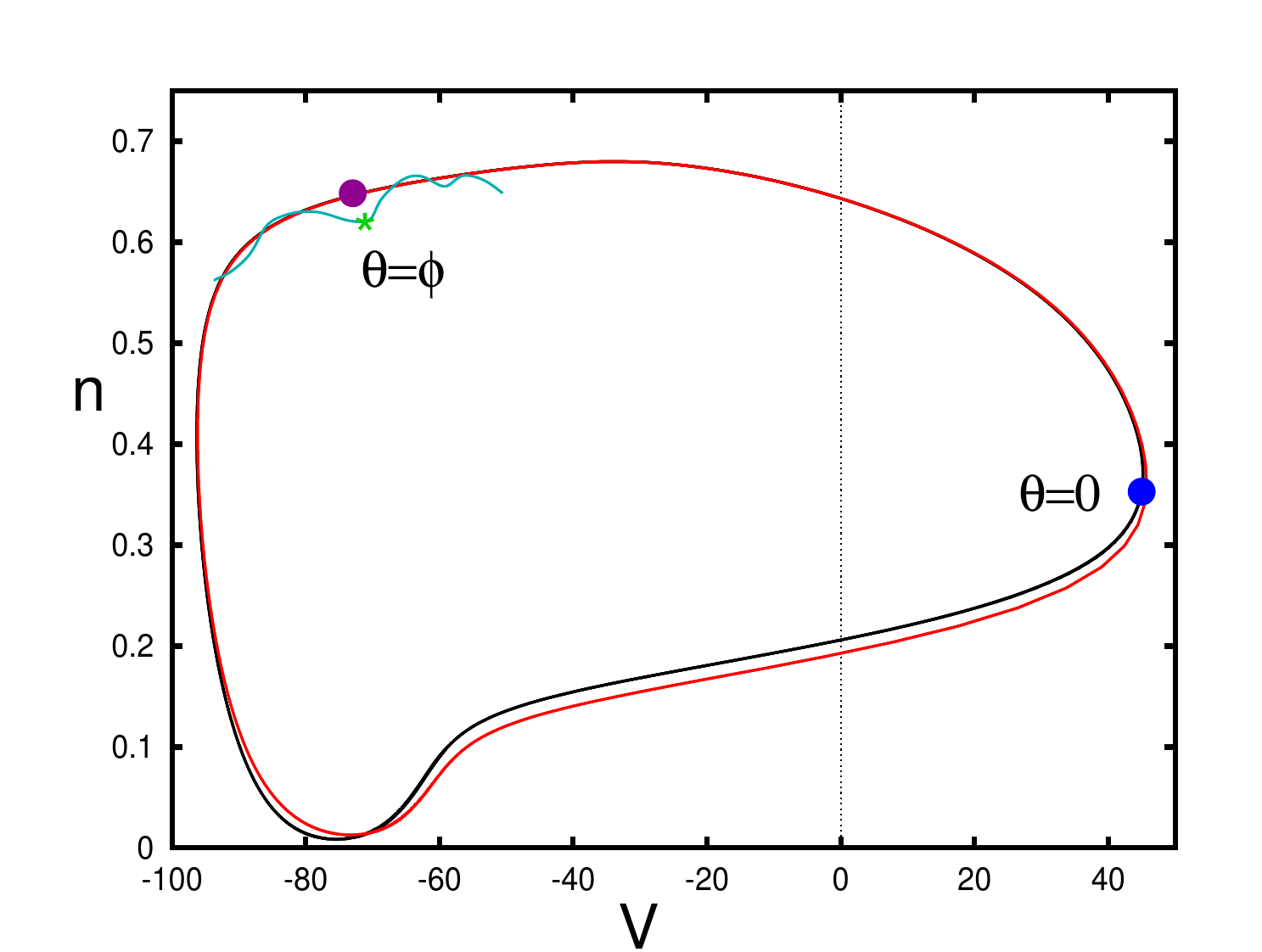}
\caption{How the phase is extracted.  A trajectory in green is shown near the limit cycle (red, for $q=0.1$). The point (asterisk) on the trajectory is closest to the magenta point on the limit cycle which has phase $\phi$, so we assign this phase.  The black curve is the projection of the limit cycle for $q=0.5$.  }
\label{fig:trbext}
\end{figure}

\ymp{In order to compare the slowly varying phase model to the full model, we need a way to extract the phase from the full model.  For the simple $\lambda-\omega$ model, we could get the exact phase since the limit cycle is circle with a constant angular velocity. One method that is commonly used is to apply a Hilbert transform to the voltage and then extract the phase from this.  However, for the Traub model (and, in fact, any model), we have more than just the voltage, so we can extract an approximate phase by picking a point on the unperturbed limit cycle that is closest to the point whose phase we wish to determine. (This is a fairly crude approximation; ideally, we would determine which isochron the point lies on by integrating the initial data forward for several periods and then matching the point. This method is very time consuming \cite{danzl}, so we have opted for the simpler approximation.) Figure \ref{fig:trbext} shows how this is done. We take the $(V,n)$ coordinates of the simulation and find the value of $(V,n)$ on the projected limit cycle that is closest in distance to the point on the actual trajectory. Since the voltage ($V$) spans a region of about 150 and the recovery ($n$) spans values between 0 and 1, we scale the distance metric accordingly.  We compute the variance of $V_0(t),n_0(t)$ over one cycle of the unperturbed limit cycle, call these $(\sigma_V^2,\sigma_n^2).$  Thus the distance is:
\[
\hbox{dist}(\Delta V,\Delta n):=\sqrt{(\Delta V)^2/\sigma_V^2+(\Delta n)^2/\sigma_n^2}.
\]     
We define the phase of a point $(V(t),n(t))$ to be the value $\phi$ that minimizes:
\[
\hbox{dist}(V(t)-V_0(\phi T/(2\pi)),n(t)-n_0(\phi T/(2\pi))),
\]
where $T$ is the natural period of the unperturbed limit cycle.  We pick the comparison limit cycle $(V_0,n_0)$ for a fixed value of the slowly varying parameter that is the mean value. However, as the figure shows, the phase portrait is very similar for two different values of $q.$ 
As we will see later, this method produces a very reasonable approximation of the phase.}

\begin{figure}[h!]
 \includegraphics[width=\textwidth]{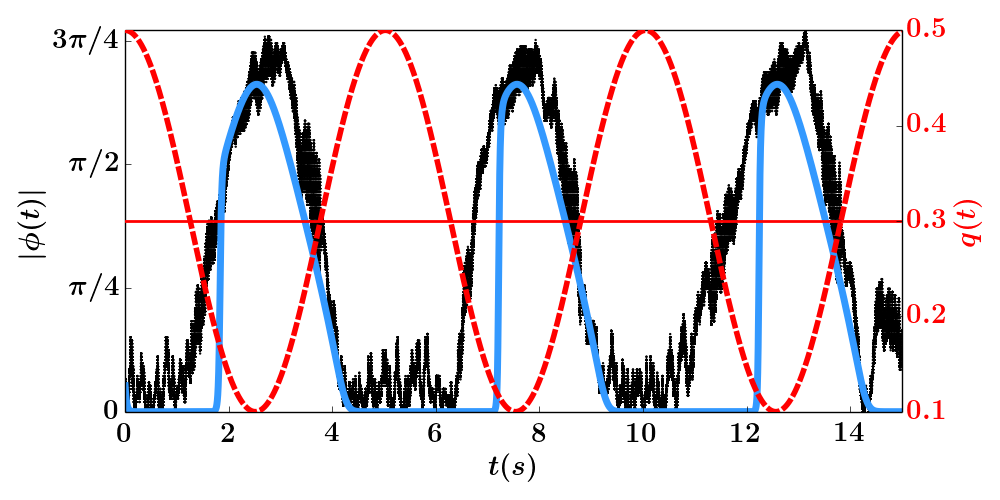}
 \caption{Periodic slowly varying parameter. Absolute value of phase difference $|\phi| = |\theta_2 - \theta_1| \ymp{\in [0,2\pi)}$ theory (light blue) vs numerics (black dots). The slow periodic parameter is shown as a dashed red line. The horizontal line represents the parameter value $q$ at which there is onset or offset of synchrony. \ymp{$\varepsilon=0.0025$, $f=5$. 245 cycles.}}\label{fig:trb_p}
\end{figure}

The next three figures compare the approximated phase model with the approximated phase extracted from the full model. That is, there are several levels of approximation to compare the theory to the full model. As described above, we approximate the function $G(\phi,\tau)$ in equation (\ref{eq:phi}) by two sine terms whose coefficients are $\tau-$dependent (cf Figure \ref{fig:trbfg}b).  We apply the same slowly varying function for the conductance of the adaptation current to get the $\tau-$dependence for the phase model. We extract the approximate phase-difference from the full model and compare the result with the phase-difference derived from equation (\ref{eq:phi}).  Figure \ref{fig:trb_p} shows the result of letting $q$ vary periodically in time; the period is 5 seconds.  The dashed red curves show the modulation and the red line shows the mean value.  The light blue curve is the phase-difference as predicted by equation (\ref{eq:phi}) and the black dots are the instantaneous approximate phase-differences from the model equations. \ymp{Each dot represnts the phase value at approximately 1/300 of one period. Because the period of the oscillation varies from 12.65s to 24.6s as a function of the slowly varying parameter, we can not give a precise total number of cycles. However, based on the total times one oscillator passes through zero phase, we estimate that there are 245 total cycles.}

\begin{figure}
\includegraphics[width=\textwidth]{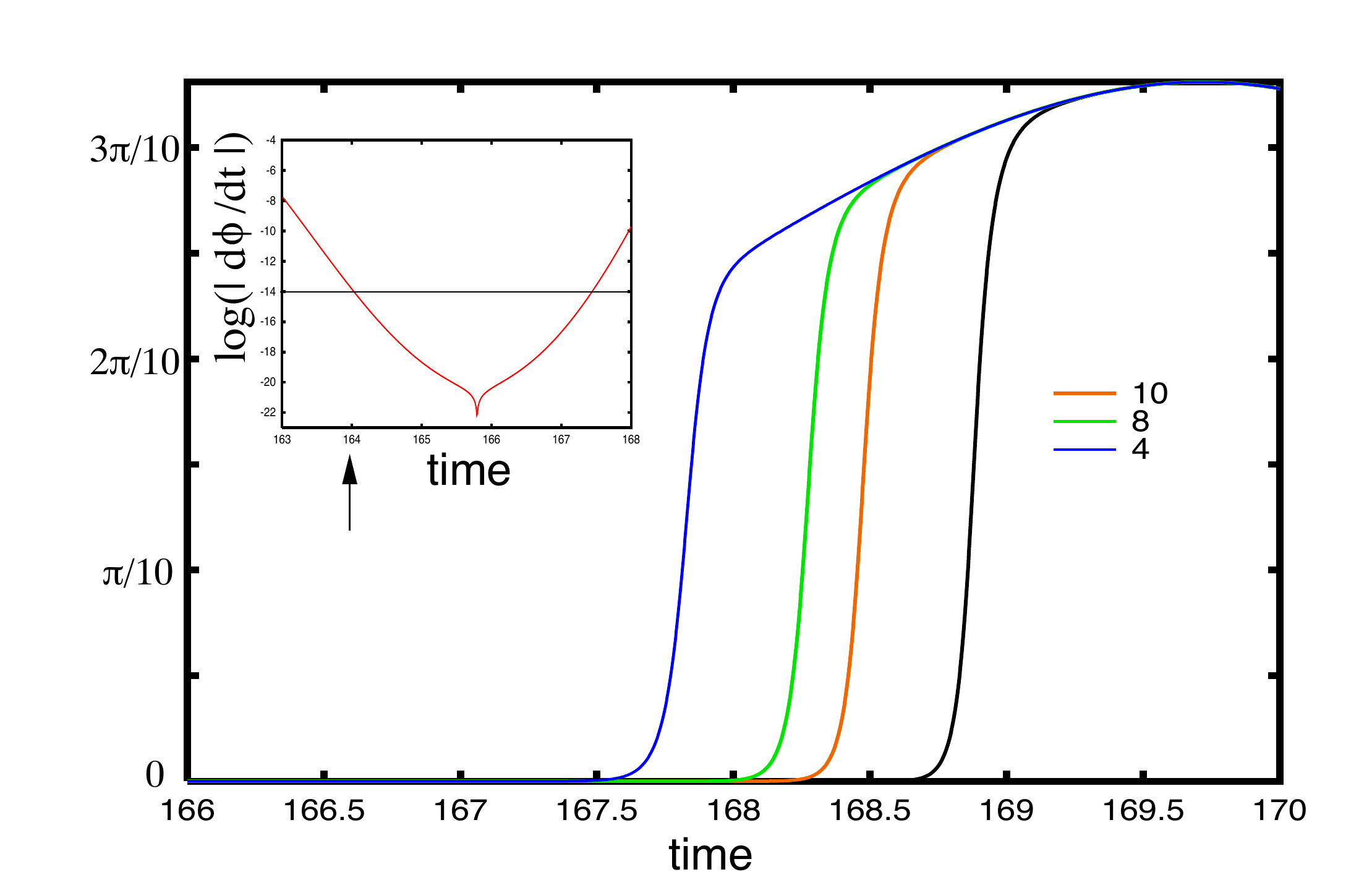}
\caption{Small perturbations of the phase model near the transition. At $t=164$, the value of $\phi$ is increased for $10^{-14}$ to $10^{-10},10^{-8},10^{-4}$ (red,green,blue), leading to an earlier jump-up time. Inset shows the log of $d\phi/dt.$ 
}\label{fig:expl}
\end{figure}

\begin{figure}[h!]
 \includegraphics[width=\textwidth]{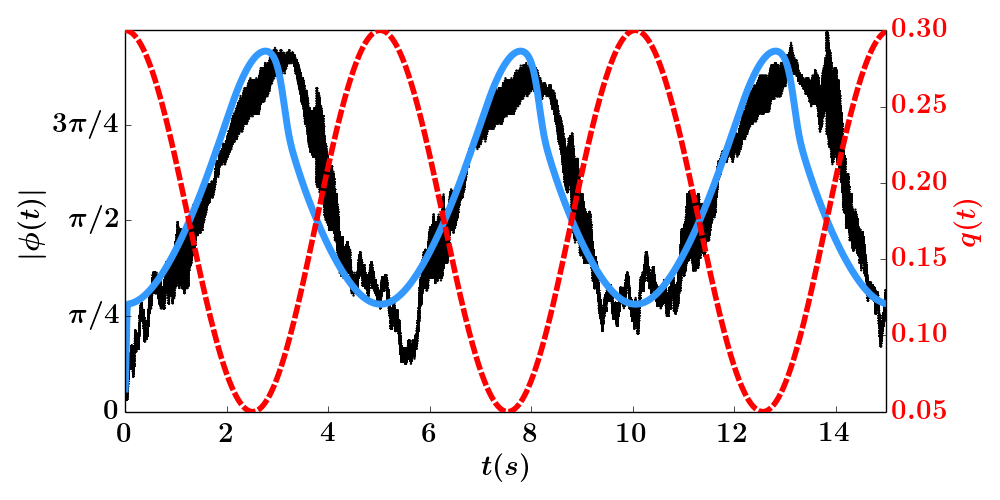}
 \caption{Periodic slowly varying parameter. Absolute value of phase difference $|\phi| = |\theta_2 - \theta_1| \ymp{\in [0,2\pi)}$ theory (light blue) vs numerics (black dots). The slow periodic parameter is shown as a dashed red line.\ymp{ The slowly varying parameter constants are $q_0=0.175$, $q_1=0.125$. $\varepsilon=0.0025$, $f=5$. 219 cycles.}}

\label{fig:trb_p_n}
\end{figure}

On the falling phase of the modulation (say, \ymp{$t=2.5-5$, $t=8-10$}, etc) the phase model and the full model agree very closely,  On the rising phase, the reduced system lags the full system by quite a bit.  Since both the rising and falling parts of the stimulus include all ranges of $q$, this difference cannot be due to a bad approximation of the interaction functions. As we noted above, the synchronous solution is a fixed point and for a range of $q$, it is attracting. Because synchrony is a fixed point and we are slowly changing from stable to unstable, there is great sensitivity at the transition.  Small changes (such as ignoring small higher order terms in the perturbation) can have drastic effects on the ``jump-up'' time as synchrony loses stability.  This is an example of a slow passage through a bifurcation \cite{maree1996slow}.   To see what we mean here, we simulate the phase-model with the periodic stimulus and perturb the phase-difference, $\phi$ by slightly increasing it when it is close to 0.  Figure \ref{fig:expl} shows the result of such a manipulation.  By increasing $\phi(t=164)$ from, say, $10^{-14}$ to $10^{-4}$ (this is still an order of magnitude smaller than the $\epsilon$ used in the simulations), we can advance the ``jump-up'' time by almost an eighth of the cycle.  The inset of the figure shows that $d\phi/dt$ is very small at this point.  For this reason, we can expect that the main error will be on the up-jump since $\phi$ has to escape from the equilibrium point at zero.  We will see similar, although less drastic, effects in the subsequent comparisons. By reducing the range of the slow parameter so that it is never close to the value for which synchrony is an attractor, we can do a much better job of tracking the phase-difference through the reduced model. Figure \ref{fig:trb_p_n} shows an example where the modulated adaptation never gets to a region where synchrony is stable. In this case, the phase-difference for the phase-reduced model never gets close to 0 and the modulation stays away from any bifurcation points.

\begin{figure}[h!]
 \includegraphics[width=\textwidth]{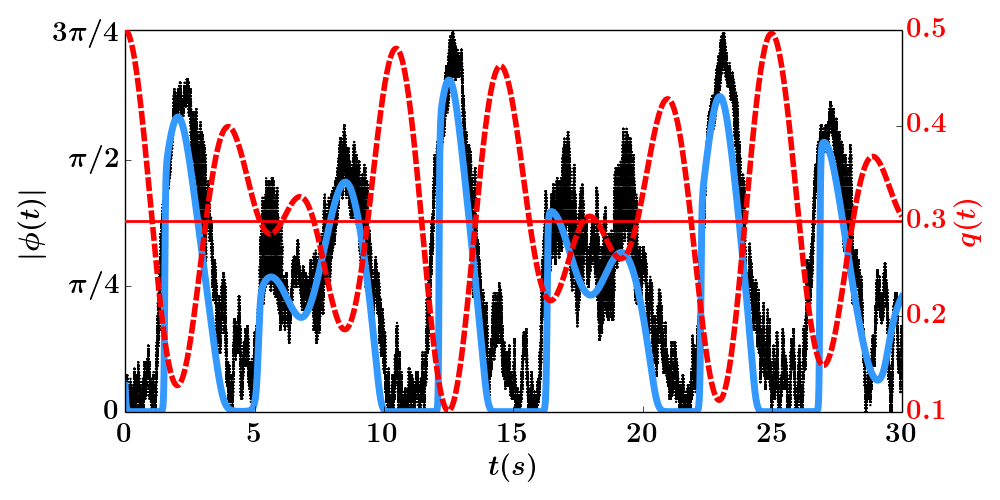}
 \caption{Quasi-periodic slowly varying parameter. The absolute value of phase difference $|\phi| = |\theta_2 - \theta_1| \ymp{\in [0,2\pi)}$ in theory (light blue) vs numerics (black dots). The quasi-periodic parameter is shown as a dashed red line. \ymp{$\varepsilon=0.0025$, $f=5$. 444 cycles.}}\label{fig:trb_qp}
\end{figure}

Figure \ref{fig:trb_qp} is similar to figure \ref{fig:trb_p}, except that the modulation is quasi-periodic. As with the periodic modulation, the phase model follows the full model quite closely once the system jumps away from the synchronous equilibrium.  However, like the periodic case, the phase model has a delayed jump-up from synchrony relative to the full model; this is especially evident at \ymp{$t\approx 25$}.

\begin{figure}[h!]
 \includegraphics[width=\textwidth]{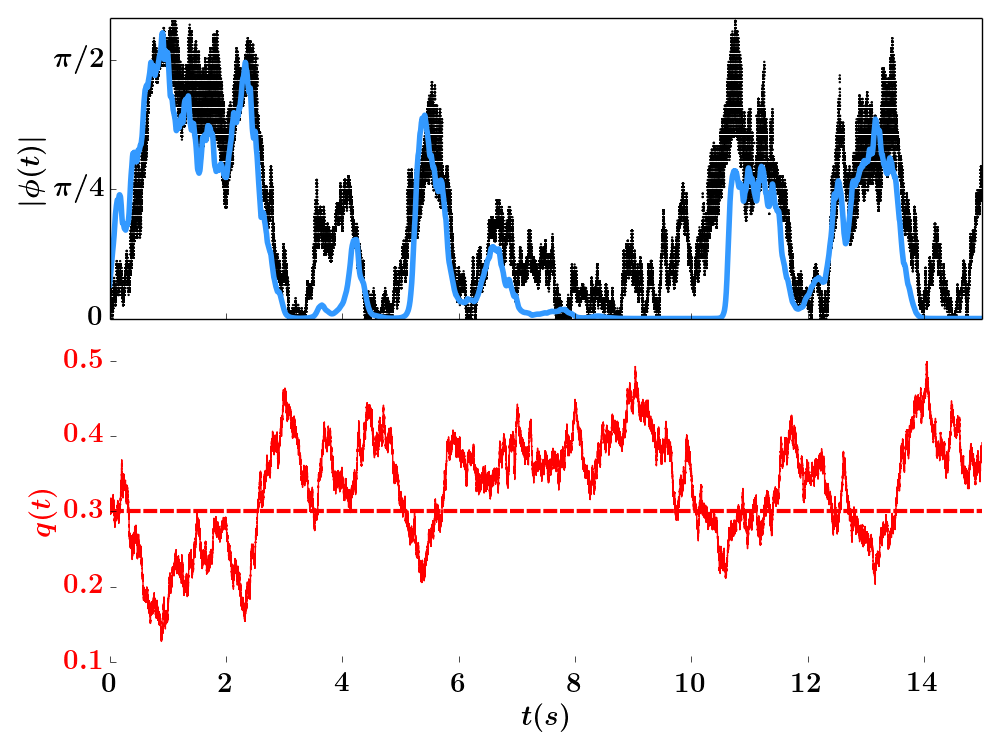}
\caption{Noisy slowly varying parameter. Absolute value of phase difference $|\phi| = |\theta_2 - \theta_1| \ymp{\in [0,2\pi)}$ theory (light blue) vs numerics (black dots). The noisy parameter is shown as a dashed red line. XPP seed 4. \ymp{$\varepsilon=0.0025$, $f=5$. 404 cycles.}}\label{fig:trb_s}
\end{figure}

Finally, in Figure \ref{fig:trb_s}, we use a slowly varying stochastic signal that is generated by an Ornstein-Uhlenbeck process and then rescaled so that the range is $[-1,1].$  As in figures \ref{fig:trb_p} and \ref{fig:trb_qp}, the phase model does a fairly good job of tracking the full model. Similarly, the jump up from synchrony is often delayed  (especially evident for \ymp{$t\in [2.5,4]$} ) as was the case in all the previous simulations.

The slowly modulated interaction function works well in spite of the many approximations that we have made in the biophysical model.

\subsection{Networks and synchrony}

\begin{figure}[h!]
 \includegraphics[width=\textwidth]{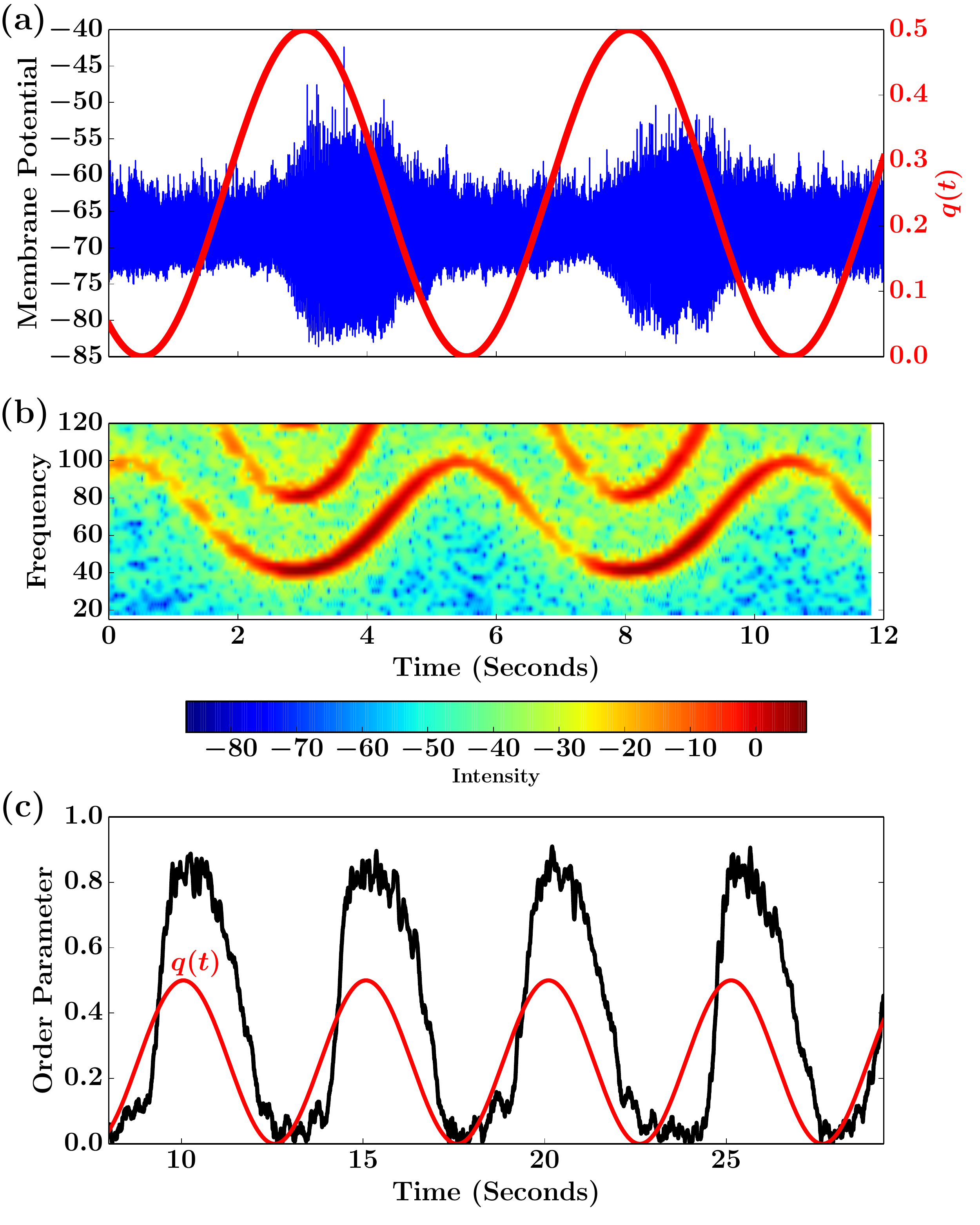}
\caption{Network of 51 Traub oscillators all-all coupled as in previous figures. (a) Summed voltage $V_{tot}=\sum_j V_j(t)$ along with the modulation of the potassium conductance; (b) Spectrogram  showing greatly increased power when the conductance is high; (c) Order parameter from the phase reduction showing a similar increase 
}
\label{fig:trb50}
\end{figure}

The methods we have described have, so far, been applied only to two coupled oscillators.  There is no reason why we cannot apply them to networks as well.  In this case, it is interesting to consider the idea of global synchronization in the presence of modulation.  Here we consider (for simplicity) a population of $N$ (we take $N=51$, here) globally coupled neurons that are subject to slow modulation of the $M-$current as in the previous sections.  We weakly couple the  Traub model neurons with excitatory coupling and slow periodic modulation of the adaptation. Coupling is all-all and divided by the total number of neurons. Thus, each includes the synaptic current, $I_{syn}=g_{syn}s_{tot}(t)(V-E_{syn})$, where
\begin{equation}
\label{eq:stot}
 s_{tot}(t)=\frac{1}{51} \sum_{j=0}^{50}(t)s_j(t)
\end{equation}
and $s_j(t)$ are the individual synaptic gating variables for each neuron.
Figure \ref{fig:trb50} shows the result of the simulation. As a surrogate for, say, the local field potential, we look at the total voltage of all the oscillators, $V_{tot}=(1/N)\sum_j V_j(t).$  Panel A shows the full picture of $V_{tot}(t)$ over 12 seconds.  It is difficult to see the synchronization, but can be roughly judged by looking at the variance of $V_{tot}$: larger variance means greater synchrony. \ymp{(If the oscillators were completely asynchronous, their sum would be close to a constant and so the variance of the sum will be small. If they are completely synchronized, then the variance of the sum will be large as the voltage swings over a 150 mV range.)} To better illustrate this point, we have also computed the spectrogram (panel B) over this period of time. Notice the large red band that starts at the peak of $q(t)$ and tails off as $q(t)$ tends to zero.  Higher bands represent harmonics of the oscillations. This panel also illustrates the dramatic effect that adaptation has on the frequency of the rhythm which ranges between 40 and 100 Hz. Higher frequencies correspond to lower adaptation and weaker synchrony.  
We can apply the same phase reduction methods to this model to get a system of phase equations:
\[
\theta_i' = \frac{1}{N+1} \sum_{j=0}^{N} H(\theta_j-\theta_i,\tau) + \sigma \xi_j
\]
where we have added some weak noise, $\sigma$ to push off the invariant synchrony manifold.  To quantify the synchronization, we look at the order parameter:
\[
\hbox{OP} = \frac{1}{N+1} \left|\sum_{j=0}^{N}e^{i\theta_j}\right|.
\]
Figure \ref{fig:trb50}c shows the clear periodic waxing and waning of OP as the slowly varying potassium conductance goes from large to small. When $q(t)$ is close to zero (no adaptation), the OP is also near zero and as \ymp{$q(t)$} tends to its maximum value of $0.5$, OP gets very close to 1.    Thus, we see that slow modulation of this type of network shows transitions in and out of synchrony. We expect similar effects for non-periodic modulation as long as it is sufficiently slow.

\section{Discussion}
We have shown that it is possible to still accurately apply weak coupling theory and phase reduction to oscillators even in a changing environment when the changes are occurring at a sufficiently slow time scale. In a previous paper \cite{rubinrubin}, we showed that slow {\em coupling} between oscillators was equivalent to weak coupling and that when the slow parameters were frozen, then the spike-to-spike synchrony in some moderate interval of time could be predicted by the corresponding interaction function. Here, we formalize this notion and demonstrate that the interaction functions are time-dependent (with respect to the slowly varying parameter) and thus, we do not need to freeze any parameters.   In the present theory, the slow forcing was exogenous and imposed on the system. In contrast, in bursting systems, the slow modulation is internally generated as the slow variable goes through various bifurcations between quiescence and oscillations.  In \cite{sherman1994}, Sherman studied two weakly coupled bursters and observed that during the spiking phase (when the fast system is periodic), spikes did not synchronize but were driven to asymmetric and out-of-phase oscillations.   The methods we have developed here require the existence of a limit cycle, so they cannot be applied globally to the spike synchronization during autonomous bursting. However, if we make the reasonable assumption that during the quiescent stage of the burst, the two cells are drawn to the slowly changing equilibrium, then, when they jump up, they will be nearly synchronous (but, not quite; if they were precisely synchronous, then they would stay that way for all time). We can then use the theory developed here to study the dynamics of the spike-to-spike synchrony during the active period of the burst and thus explain in a more formal matter Sherman's results.      

Gutkin et al \cite{gutkinreyes} looked at how spiking neural oscillators subjected to a slowly varying input responded to brief perturbations given between spikes. They measured the inter-spike interval (ISI) without the perturbation and then fixed the input so the oscillator had the same period as the ISI. With this fixed oscillator, they computed the phase response curve (the function $Z(t)$ solving the adjoint problem) and used this to predict how the perturbation would affect the spike time in the slowly varying system.  The method developed here, could be used to improve this estimate since we know both the slowly varying frequency and the slowly varying function $Z(t;\tau)$.  This type of correction was, in fact, the goal of \cite{kureb13}.   

Slowly varying inputs differ in ways both quantitative and, more importantly, qualitative from faster inputs.   For example, suppose that two oscillators receive identical periodic inputs that have a frequency close to the unforced frequency of the oscillators.  Then for some range of input amplitudes, we can expect the oscillators to lock in a 1:1 manner with each other and thus be completely synchronized even in absence of coupling.  Similarly, weak identical noise applied to two uncoupled oscillators will also synchronize them \cite{pikovsky,teramae2004,ermentrout2008}, but the noise has to be sufficiently fast; synchrony falls off rapidly as the time constant of the noise slows down \cite{galan2008}.  Thus, fast common rapidly changing inputs will tend to synchronize uncoupled oscillators. But the slow modulations we study here have no such properties. Indeed, looking at equations (\ref{eq:tha}-\ref{eq:thb}), the common {\em slow} input cannot move the phase-difference without direct coupling.  It would be interesting to look at the synchronization between two slowly varying oscillators that are subjected to fast correlated noise and derive some equations for the expected phase-difference.

\section{Conclusion}
The Fredholm alternative provides a useful proof method to re-derive the phase equation in Kurebayashi et al. After obtaining the phase equation, we use the theory of weakly coupled oscillators to derive the interaction function, from which we can study the stability of synchrony and anti-synchrony.

Despite the phase estimation and the mode truncation, our theory accurately predicts the phase of the Traub model with periodic, quasi-periodic, and stochastic slowly varying parameters (we have similar positive results for the $\lambda-\omega$ system). Because the mode truncation depends on the accuracy of the numerically derived interaction functions, and because the interaction functions in turn depend only on the coupling terms and the iPRC, we can apply the mode truncation method (and subsequently our result) to any autonomous system for which the iPRC and coupling terms are known. The methods here show that we can extend the notion of weak coupling and synchronization of nonlinear neural oscillators to the more realistic scenario in which the environment is changing.

\bigskip
\noindent {\bf Acknowledgments} 
BE and YMP were partially supported by NSF DMS 1219753
\bigskip

\section{Traub Model With Adaptation}
All other equations for the Traub model are defined as follows
\begin{equation*}
\begin{split}
t_w(V)&=\tau_w/(3.3\exp((V-V_{wt})/20)+\exp(-(V-V_{wt})/20))\\
w_\infty(V)&=1/(1+\exp(-(V-V_{wt})/10))\\
a_m(V)&=0.32(54+V)/(1-\exp(-(V+54)/4))\\
b_m(V)&=0.28(V+27)/(\exp((V+27)/5)-1)\\
a_h(V)&=0.128\exp(-(V-V_{hn})/18)\\
b_h(V)&=4/(1+\exp(-(V+27)/5))\\
a_n(V)&=0.032(V+52)/(1-\exp(-(V+52)/5))\\
b_n(V)&=0.5\exp(-(57+V)/40)\\
\alpha(V)&=a_0/(1+\exp(-(V-V_t)/V_s))
\end{split}
\end{equation*}
\begin{table}
\caption {Traub parameter values} \label{tab:traubparms}
\begin{center}
\begin{tabular}{l|l}
Parameter & Value\\
\hline
$C$& $1 \mu \text{F}/\text{cm}^2$\\
$g$& $5 \text{mS}/\text{cm}^2$\\
$\varepsilon$&0.0025\\
$f$& $0.5 (\frac{1000 \varepsilon}{2\pi} \text{Hz})$\\
$I$& $3 \mu \text{A}/\text{cm}^2$\\
$V_{wt}$&$-35 \text{mV}$\\
$\tau_w$&$100 \text{ms}$\\
$E_k$&$-100 \text{mV}$\\
$E_{Na}$& $50 \text{mV}$\\
$E_l$& $-67 \text{mV}$\\
$g_l$& $0.2 \text{mS}/\text{cm}^2$\\
$g_k$& $80 \text{mS}/\text{cm}^2$\\
$g_{Na}$& $100 \text{mS}/\text{cm}^2$\\
$V_{hn}$& $-50 \text{mV}$\\
$a_0$&4\\
$\tau$& $4 \text{ms}$\\
$V_t$& $0 \text{mV}$\\
$V_s$& $5 \text{mV}$\\
$E_{syn}$& $0 \text{mV}$\\
$q_0$& $0.3 \text{mS}/\text{cm}^2$\\
$q_1$& $0.2 \text{mS}/\text{cm}^2$\\
\end{tabular}
\end{center}
\end{table}
\subsection{Fourier Coefficients}
\label{app:fourier_coeff}
The Fourier coefficients used in the approximation are shown in Table~\ref{tab:trb_fourier}
%
% \begin{table}[h!]
% \caption{Traub Fourier Sine Coefficients} \label{tab:trb_s_fourier}
% \begin{center}
% \begin{tabular}{l|l}
%  Parameter & Value\\
%  \hline
%  $b_1(0.1)$ & 0.721387113706\\
%  $b_1(0.3)$ & -1.5028098729\\
%  $b_2(0.1)$ & 0.738312597998\\
%  $b_2(0.3)$ & 1.03494013487
% \end{tabular}
% \end{center}
% \end{table}
% %

\begin{table}[h!]
\caption{Traub Fourier Coefficients} \label{tab:trb_fourier}
\begin{center}
\begin{tabular}{l|l}
 Cosine & Sine\\
 \hline
  $a_0(0.1)= 19.6011939665$ & -\\
  $a_0(0.3)= 17.4255017198$& -\\
  $a_1(0.1)=-3.32476526025$ & $b_1(0.1)=0.721387113706$\\
  $a_1(0.3)=-6.97305767558$ & $b_1(0.3)=-1.5028098729$\\
  $a_2(0.1)=-0.255371105623$ & $b_2(0.1)=0.738312597998$\\
  $a_2(0.3)-0.83690237427$ & $b_2(0.3)=1.03494013487$
\end{tabular}
\end{center}

\end{table}

\bibliographystyle{plain}
\bibliography{bibliography}

\end{document}